\documentclass{conm-p-l}
\usepackage{amssymb,colordvi}
\usepackage{epsf}

\copyrightinfo{2003}{(American Mathematical Society)}

\numberwithin{equation}{section}

\newtheorem{theorem}[equation]{Theorem}

\newtheorem{cor}[equation]{Corollary}

\theoremstyle{definition}
\newtheorem{definition}[equation]{Definition}

\newtheorem{example}[equation]{Example}

\theoremstyle{remark}
\newtheorem{remark}[equation]{Remark}
\newcounter{FNC}[page]
\def\newfootnote#1{{\addtocounter{FNC}{2}$^\fnsymbol{FNC}$%
     \let\thefootnote\relax\footnotetext{$^\fnsymbol{FNC}$#1}}}

\begin{document}

\title[Real toric varieties, and the algebraic moment map]%
      {Toric ideals, real toric varieties,\\ and the algebraic moment map}

\author{Frank Sottile}
\address{Department of Mathematics\\
         Texas A\&M University\\
         College Station\\
         TX \ 77843\\
         USA}
\email{sottile@math.tamu.edu}
\urladdr{www.math.tamu.edu/\~{}sottile}
\thanks{Research supported in part by NSF grant DMS-0134860}
\thanks{\Red{Corrected version of published article}}

\subjclass[2000]{14M25, 14Q99, 13P10, 65D17, 68U05, 68U07}
\thanks{1998 ACM Computing Classification System:  
           I.3.5 Computational Geometry and Object Modeling}
\keywords{Toric Varieties, Moment Map, B\'ezier Patches, Linear Precision}

\begin{abstract}
 This is a tutorial on some aspects of toric varieties related to their
 potential use in geometric modeling.
 We discuss projective toric varieties and their ideals, as well as real
 toric varieties.
 In particular, we explain the relation between linear precision and a
 particular linear projection we call the algebraic moment map. 
\end{abstract}

\maketitle
\section*{Introduction}

We develop further aspects of toric varieties that may be useful in geometric
modeling, building on Cox's introduction to toric varieties, {\sl What is a
toric variety?}~\cite{Co02}, which also appears in this volume.
Notation and terminology follow that article, with a few small exceptions.
This paper is organized into eight sections:

\begin{enumerate}
 \item[1.] Projective Toric Varieties 
 \item[2.] Toric Ideals
 \item[3.] Linear Projections
 \item[4.] Rational Varieties 
 \item[5.] Implicit Degree of a Toric Variety
 \item[6.] The Real Part of a Toric Variety
 \item[7.] The Double Pillow
 \item[8.] Linear Precision and the Algebraic Moment Map
\end{enumerate}

\section{Projective Toric Varieties}
In this tutorial, we study toric varieties as subvarieties of projective
space.
This differs slightly from Cox's~\cite{Co02} presentation, where 
toric varieties are studied  via the abstract toric variety
$X_\Sigma$ of a fan\index{fan} $\Sigma$.
The resulting loss of generality is compensated by  
the additional perspective this alternative view provides. 
Only in the last few sections do we discuss abstract toric varieties.

A projective toric variety\index{toric variety} may be given as the closure of
the image of a map 
 \[
   (\mathbb{C}^*)^n\ \longrightarrow\ \mathbb{P}^\ell\,,
 \]
defined by Laurent monomials as in Section 13 of~\cite{Co02}.
There, the monomials had exponent vectors given by all the integer lattice
points in a polytope. 
Here, we study maps given by {\it any} set of Laurent monomials.

Our basic data structure will be a list of integer exponent vectors
 \[
   \mathcal{A}\ :=\ \{ {\bf m}_0, {\bf m}_1,\dotsc,{\bf m}_\ell\}\ 
                     \subset\ \mathbb{Z}^n\,.
 \]
Such a list gives rise to a map 
$\varphi_{\mathcal{A}}$ (written $\varphi$ when $\mathcal{A}$ is understood), 
 \begin{equation}\label{E:calA}
  \begin{array}{rcccl}
    \varphi_{\mathcal{A}}&\colon&
                 (\mathbb{C}^*)^n&\longrightarrow&\mathbb{P}^\ell\\
    &&{\bf t}&\longmapsto&[{\bf t}^{{\bf m}_0}, 
                         {\bf t}^{{\bf m}_1},\dotsc,{\bf t}^{{\bf m}_\ell}]\,.
  \end{array}
 \end{equation}
We explain this notation.
Given ${\bf t}=(t_1,\dotsc,t_n)\in(\mathbb{C}^*)^n$ and an exponent vector
${\bf a}=(a_1,\dotsc,a_n)$, the monomial 
${\bf t}^{\bf a}$ is equal to $t_1^{a_1}t_2^{a_2}\dotsb t_n^{a_n}$.
In this way, the coordinates of $\mathbb{P}^\ell$ are naturally indexed by the
exponent vectors lying in $\mathcal{A}$.
The {\it toric variety} $Y_{\mathcal{A}}$ is the closure in $\mathbb{P}^\ell$ of
the image of the map $\varphi_{\mathcal{A}}$.
This map $\varphi_{\mathcal{A}}$ gives a parametrization of $Y_{\mathcal{A}}$ 
by the monomials whose exponents lie in $\mathcal{A}$.

We claim that $Y_{\mathcal{A}}$ is a toric variety as defined in Section 2
of~\cite{Co02}.
The map 
\[
   (\mathbb{C}^*)^n\ \ni\ {\bf t}\ \ \longmapsto
   \ \ ({\bf t}^{{\bf m}_0},{\bf t}^{{\bf m}_1},\dotsc,{\bf t}^{{\bf m}_\ell})\ 
    \in\ (\mathbb{C}^*)^{1+\ell}
\]
is a homomorphism from the group $(\mathbb{C}^*)^n$ to the group
$(\mathbb{C}^*)^{1+\ell}$ of  diagonal $(1+\ell)$ by $(1+\ell)$ matrices,
which acts on $\mathbb{P}^\ell$.
Thus $(\mathbb{C}^*)^n$ acts on $\mathbb{P}^\ell$ via this map.
Since scalar matrices (those in $\mathbb{C}^*I_{1+\ell}$) act trivially on 
$\mathbb{P}^\ell$, this action of $(\mathbb{C}^*)^{1+\ell}$  on 
$\mathbb{P}^\ell$ factors through the group 
$(\mathbb{C}^*)^{1+\ell}/\mathbb{C}^*I_{1+\ell}\simeq(\mathbb{C}^*)^\ell$, 
which is the dense torus in the toric variety $\mathbb{P}^\ell$.
Then $Y_{\mathcal{A}}$ is the closure of the image of $(\mathbb{C}^*)^n$ in
this torus, that image $T$ acts on $Y_{\mathcal{A}}$, and thus $T$
is the dense torus of $Y_{\mathcal{A}}$.\smallskip

Suppose that $\mathcal{A}=\Delta\cap\mathbb{Z}^n$, where $\Delta$ is a lattice
polytope. 
Then $Y_{\mathcal{A}}$ is the image of the abstract toric variety
$X_\Delta$ given by the normal fan\index{fan!normal, of a polytope} of $\Delta$
under the map of Section 13 in~\cite{Co02}. 
When $\mathcal{A}$ has this form, we write $Y_\Delta$ for $Y_{\mathcal{A}}$.

\begin{example}\label{E:3examples}
 Consider the three lattice polytopes
 \begin{align*}
  [n],\   &\textrm{the line segment } [0,n]\subset \mathbb{R}\,,\\
  \triangle_n, \ 
   &\textrm{the triangle } \{(x,y)\in\mathbb{R}^2\mid 0\leq x,y,\ 
                                 x+y\leq n\}\,,\quad \textrm{and}\\
  \Box_{m,n},\ &\textrm{the rectangle }
            \{(x,y)\in\mathbb{R}^2\mid 0\leq x\leq m,\ 
                                 0\leq y\leq n\}\, .
 \end{align*}
 The maps $\varphi$ for these polytopes are
 \begin{eqnarray*}
    t\hspace{8.1pt}
     &\longmapsto&[1,t,t^2,\dotsc,t^n]\ \in\ \mathbb{P}^n\,,\\
   (s,t)&\longmapsto&[1,s,t,s^2,st,t^2,\dotsc,s^n,s^{n-1}t,\dotsc,t^n]\ 
         \in\ \mathbb{P}^{\binom{n+2}{2}}\rule{0pt}{11pt}\,,\quad \textrm{and}\\
    (s,t)&\longmapsto&
      [1,s,\dotsc,s^m,t,st,\dotsc,s^mt,\dotsc,t^n,st^n,s^mt^n]\ \in\ 
           \mathbb{P}^{mn}\rule{0pt}{11pt}\, .
 \end{eqnarray*}
and the resulting projective toric varieties are known (see~\cite{Ha92}) as
 \begin{eqnarray*}
  Y_{[n]} &=&\textrm{the rational normal curve in }\mathbb{P}^n\,,\\
  Y_{\triangle_n} &=&\textrm{the Veronese embedding of $\mathbb{P}^2$ in 
              $\mathbb{P}^{\binom{n+2}{2}}$}\,,\quad \textrm{and}\\
  Y_{\Box_{m,n}}&=&\textrm{the Segre embedding of 
                     $\mathbb{P}^1\times\mathbb{P}^1$ in 
               $\mathbb{P}^{mn}$ of bidegree $m,n$}\, .
 \end{eqnarray*}
 In geometric modeling these projective toric varieties give rise to,
 respectively, B\'ezier curves\index{B\'ezier curve}, rational B\'ezier
 triangles of degree $n$, and\index{B\'ezier patch} 
 tensor product surfaces of bidegree $(m,n)$.
 \hfill \raisebox{-3pt}{\epsfxsize=11pt\epsffile{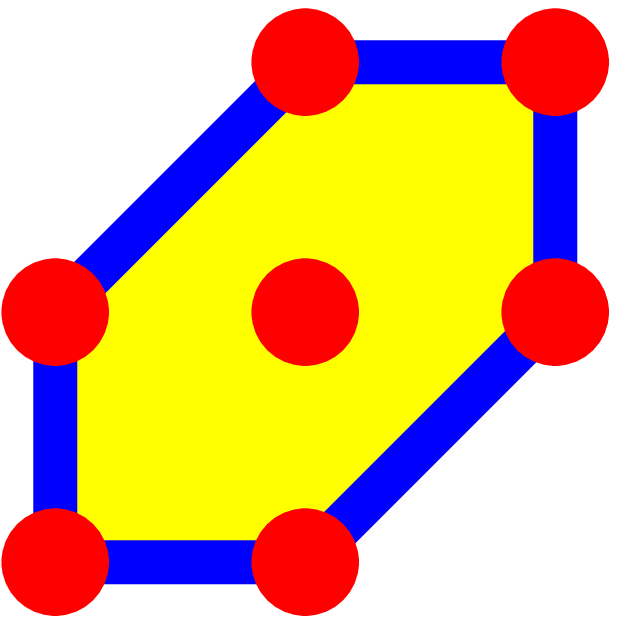}}
\end{example}

\begin{example}
 Let $n=1$ and $\mathcal{A}=\{0,2,3\}$.
 Then the map~\eqref{E:calA} is
\[
   t\ \longmapsto\ [1, t^2, t^3]
\]
 whose image $Y_{\mathcal{A}}$ is the cuspidal cubic
\[
  \raisebox{0.392in}{$ {\bf V}(x_0x_2^2-x_1^3)\ \ =\ \ \ \ $} 
   \epsfysize=0.85in\epsffile{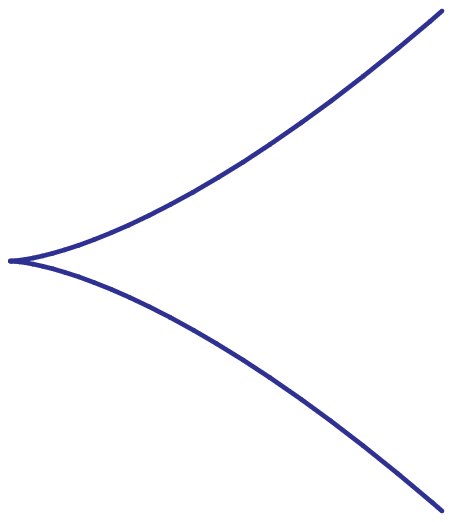}  
\]
 which is the non-normal toric variety of Example 3.2 in~\cite{Co02}.
 \hfill \raisebox{-3pt}{\epsfxsize=11pt\epsffile{figures/QED.eps}}
\end{example}

\begin{example}\label{E:hex}
 Let $\Delta$ be the hexagon which is the convex hull of the six column vectors
 $\genfrac{[}{]}{0pt}{1}{1}{0},\genfrac{[}{]}{0pt}{1}{1}{1},
  \genfrac{[}{]}{0pt}{1}{0}{1},\genfrac{[}{]}{0pt}{1}{-1}{\ 0},
  \genfrac{[}{]}{0pt}{1}{-1}{-1},\genfrac{[}{]}{0pt}{1}{\ 0}{-1}$. 
 We depict $\Delta$ and its normal fan $\Sigma_\Delta$.
 \[
  \begin{picture}(95,65)
    \put(0,30){$\Delta$ =}
    \put(30,0){\epsffile{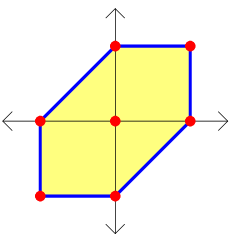}}
  \end{picture}
  \qquad
  \begin{picture}(22,65)
   \put(0,30){and}
  \end{picture}
  \qquad
  \begin{picture}(100,65)
    \put(0,30){$\Sigma_\Delta$ =}
     \put(35,0){\epsffile{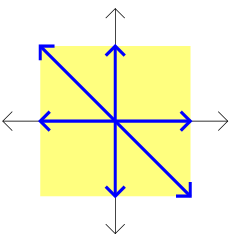}}
  \end{picture}
 \]
 Then $\Delta\cap\mathbb{Z}^2$ consists of these six vectors
 (the vertices of $\Delta$) together with the origin
 $\genfrac{[}{]}{0pt}{1}{0}{0}$.
 Thus $Y_\Delta =
   \overline{\{ [1,t,st,s,t^{-1},s^{-1}t^{-1},s^{-1}]\mid
   s,t\in\mathbb{C}^\times\}} \subset \mathbb{P}^6$.
  \hfill \raisebox{-3pt}{\epsfxsize=11pt\epsffile{figures/QED.eps}}
\end{example}

\begin{remark}\label{R:U_0}
 Suppose that the origin ${\bf 0}$ is an element of $\mathcal{A}$ and 
 that $m_0={\bf 0}$.
 Then the image of the map $\varphi$ of~\eqref{E:calA} lies in the
 {\it principal affine part} of $\mathbb{P}^\ell$
\[
  U_0\ :=\ \{ x\in\mathbb{P}^\ell\mid x_0\neq 0\}\ \simeq\ \mathbb{C}^\ell\,,
\]
 whose coordinates are $[1,x_1,x_2,\dotsc,x_\ell]$.
 Thus $U_0\cap Y_{\mathcal{A}}$ is an affine toric variety.
 In this case, the dimension of the projective toric variety $Y_{\mathcal{A}}$
 is equal to the dimension of the linear span of the exponent vectors
 $\mathcal{A}$. 
 \hfill \raisebox{-3pt}{\epsfxsize=11pt\epsffile{figures/QED.eps}}
\end{remark}

\section{Toric Ideals}
The {\it toric ideal}\index{toric ideal} $I_{\mathcal{A}}$ is the homogeneous
ideal of polynomials 
whose vanishing defines the projective toric variety
$Y_{\mathcal{A}}\subset\mathbb{P}^\ell$.
Equivalently, $I_{\mathcal{A}}$ is the ideal of all the homogeneous polynomials
vanishing on $\varphi_{\mathcal{A}}((\mathbb{C}^*)^n)$.
Our description of $I_{\mathcal{A}}$ follows the presentation in Sturmfels's book, 
\emph{Gr\"obner bases and convex polytopes}~\cite{Sturmfels_GBCP}.

Let $[x_0,x_1,\dotsc,x_\ell]$ be homogeneous coordinates for
$\mathbb{P}^\ell$ with $x_j$ corresponding to the monomial 
${\bf t}^{{\bf m}_j}$ in the map $\varphi_{\mathcal{A}}$~\eqref{E:calA}, where
$\mathcal{A}=\{{\bf m}_0,{\bf m}_1,\dotsc,{\bf m}_\ell\}$.
A monomial ${\bf x}^{\bf u}$ in these coordinates has an exponent vector 
${\bf u}\in\mathbb{N}^{1+\ell}$.
Restricting the monomial ${\bf x}^{\bf u}$  to
$\varphi_{\mathcal{A}}(t_1,\dotsc,t_n)=
  [{\bf t}^{{\bf m}_0},{\bf t}^{{\bf m}_2},\dotsc,{\bf t}^{{\bf m}_\ell}]$
yields the monomial
 \[
   {\bf t}^{u_0{\bf m}_0+u_1{\bf m}_1+\dotsb+u_\ell{\bf m}_\ell}\ .
 \]
This exponent vector is $\mathcal{A}{\bf u}$, where we consider $\mathcal{A}$ to 
be the 
matrix whose columns are the exponent vectors in $\mathcal{A}$. 
For the hexagon of Example~\ref{E:hex}, this 
is
 \[
   \left(\begin{array}{ccccrrr}0&1&1&0&-1&-1&0\\
                     0&0&1&1&0&-1&-1\end{array}\right)\ 
 \]

This discussion shows that a homogeneous binomial 
${\bf x}^{\bf u}-{\bf x}^{\bf v}$ with 
$\mathcal{A}{\bf u}=\mathcal{A}{\bf v}$ vanishes on 
$\varphi_{\mathcal{A}}((\mathbb{C}^*)^n)$ and hence lies in the toric ideal
$I_\mathcal{A}$. 
In fact, the toric ideal $I_\mathcal{A}$ is the linear span of these
binomials.

\begin{theorem}\label{T:toric_ideal}
 The toric ideal $I_\mathcal{A}$ is the linear span of all homogeneous binomials 
 ${\bf x}^{\bf u}-{\bf x}^{\bf v}$ with 
 $\mathcal{A}{\bf u}=\mathcal{A}{\bf v}$.
\end{theorem}

We obtain a more elegant description of $I_{\mathcal{A}}$ if 
the row space of the matrix $\mathcal{A}$ contains the vector
$(1,\dotsc,1)$, for then the homogeneity of the binomial 
${\bf x}^{\bf u}-{\bf x}^{\bf v}$ follows from 
$\mathcal{A}{\bf u}=\mathcal{A}{\bf v}$.
It is often useful to force this condition as follows. 

Given a list $\mathcal{A}$ of exponent vectors in $\mathbb{Z}^n$,
the {\it lift} $\mathcal{A}^+$ of $\mathcal{A}$ to $1\times\mathbb{Z}^n$ is obtained 
be prepending a component of 1 to each vector in $\mathcal{A}$.
That is, 
\[
   \mathcal{A}^+\ :=\ \{(1,{\bf m})\mid {\bf m}\in\mathcal{A}\}
\]
The matrix $\mathcal{A}^+$ is obtained from the matrix $\mathcal{A}$ by
adding a new top row of 1s.
For the hexagon, this is 
 \[
   \mathcal{A}^+\ =\ \left(\begin{array}{ccccrrr}1&1&1&1&1&1&1\\
                     0&1&1&0&-1&-1&0\\
                     0&0&1&1&0&-1&-1\end{array}\right)\ .
 \]
Here is the lifted hexagon (shaded)
\[
  \epsffile{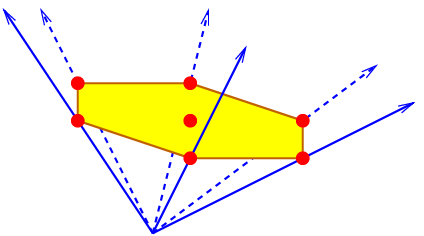}
\]
Then  $Y_{\mathcal{A}^+} = 
     \overline{\{ [r,rt,rst,rs,rt^{-1},rs^{-1}t^{-1},rs^{-1}]\mid
     r,s,t\in\mathbb{C}^\times\}} = Y_\Delta$.\smallskip 

This lifting does not alter the projective toric
variety.
Indeed, 
 \begin{eqnarray*}
  \varphi_{\mathcal{A}^+}(t_0,t_1,\dotsc,t_n)&=&
   [t_0{\bf t}^{{\bf m}_0}, t_0{\bf t}^{{\bf m}_1},\dotsc, 
       t_0{\bf t}^{{\bf m}_\ell} ]\\
 &=&
   [{\bf t}^{{\bf m}_0}, {\bf t}^{{\bf m}_1},\dotsc,{\bf t}^{{\bf m}_\ell}]\\
 &=&
   \varphi_{\mathcal{A}}(t_1,\dotsc,t_n)\,,
 \end{eqnarray*}
and so $Y_{\mathcal{A}} = Y_{\mathcal{A}^+}$.
The dimension of $Y_{\mathcal{A}}$ is one less than the dimension of 
the linear span of $\mathcal{A}^+$.
Since 
$I_{\mathcal{A}}=I_{\mathcal{A}^+}$, we have 

\begin{theorem}
 The toric ideal $I_\mathcal{A}$ is the linear span of all binomials 
 ${\bf x}^{\bf u}-{\bf x}^{\bf v}$ with 
 $\mathcal{A}^+{\bf u}=\mathcal{A}^+{\bf v}$.
\end{theorem}

If ${\bf u}\in\mathbb{Z}^{1+\ell}$, then we may write ${\bf u}$ uniquely as 
${\bf u}={\bf u}^+-{\bf u}^-$, where 
${\bf u}^+,{\bf u}^-\in{\mathbb N}^{1+\ell}$ but ${\bf u}^+$ and ${\bf u}^-$
have no non-zero components in common.
For example, if ${\bf u}=(1,3,2,-2,2,-4)$, then 
${\bf u}^+=(1,3,2,0,2,0)$ and ${\bf u}^-=(0,0,0,2,0,4)$.
We describe a smaller set of binomials that generate $I_{\mathcal{A}}$.
Let $\mbox{ker}(\mathcal{A})\subset\mathbb{Z}^{1+\ell}$ be the kernel of the
matrix $\mathcal{A}$.  

\begin{cor}\label{C:gensI}
 ${\displaystyle I_\mathcal{A}\ =\ \langle x^{{\bf u}^+}-x^{{\bf u}^-} \mid 
       {\bf u}\in\mbox{ker}(\mathcal{A}^+)}\rangle$.
\end{cor}

Algorithms for computing toric ideals are implemented in the computer algebra
systems Macaulay 2~\cite{M2} and Singular~\cite{SINGULAR}.
There are no simple formulas for a finite set of generators of a general 
toric ideal.

On the other hand, {\it quadratic binomials} in a toric ideal do have a simple 
geometric interpretation.  
Suppose that we have a relation of the form
 \begin{equation}\label{E:relation}
  {\bf a} + {\bf b}\ =\ {\bf c} + {\bf d}, \qquad\mbox{for }
   {\bf a}, {\bf b}, {\bf c}, {\bf d}\in\Delta\,.
 \end{equation}
Such a relation comes from coincident
midpoints of two line segments between lattice points in $\Delta$.
A relation~\eqref{E:relation} gives a vector 
${\bf u}\in\mbox{\it ker}(\mathcal{A})$ whose
entries are 1 in the coordinates corresponding to ${\bf a}$ and ${\bf b}$ and
$-1$ in the coordinates corresponding to ${\bf c}$ and ${\bf d}$.
The corresponding generator of the ideal $I_{\mathcal{A}}$ is $ab-cd$, where $a$
is the variable corresponding to the vector ${\bf a}$, $b$ the variable
corresponding to ${\bf b}$, and etc.

Often these simple relations suffice.
When $n=2$, Koelman~\cite{Ko93} showed that the ideal of a toric
surface $Y_\Delta$ is generated by such quadratic binomials if the polygon
$\Delta$ has more than 3 lattice points on its boundary.
Also, if a lattice polytope $\Delta\subset\mathbb{R}^n$ has the form
$n\Delta'$, for a (smaller) lattice polytope $\Delta'$, then 
the toric ideal $I_\Delta$ is generated by such quadratic binomials for
$\Delta$~\cite{BGT97}. 

\begin{example}\label{E:hex_ideal}
 Consider the ideal $I_\Delta$ of the toric variety $Y_\Delta$, where $\Delta$ is 
 the hexagon of Example~\ref{E:hex}.
 Label the exponent vectors in $\Delta$ as indicated below.
\[
 \begin{picture}(90,90)
   \put( 0, 0){\epsfxsize=90pt\epsffile{figures/hexagon.eps}}
                      \put(34,76){{\bf d}} \put(77,79){{\bf c}}
   \put( 9,52){{\bf e}} \put(49,50){{\bf a}} \put(77,50){{\bf b}}
   \put( 4,12){{\bf f}} \put(50,12){{\bf g}}
 \end{picture}
\]
There are 12 relations of the form~\eqref{E:relation}:
6 involving the midpoint of the segment connecting ${\bf a}$ and one of the
remaining vectors and 6 for the 3 antipodal pairs of points, which
all share the midpoint ${\bf a}$.
Translating these into quadratic relations gives 12 quadratic binomials in 
the toric ideal $I_\Delta$
 \begin{eqnarray*}
   &ab - cg,\ ac - bd,\  ad - ce, \ ae - df,\  af - ge, \  ag - bf,\,&\\
   &a^2 - be, \  a^2 - gd, \  a^2 - cf, \ be - cf, \ be - gd, \ cf - gd\,.&
 \end{eqnarray*}
 By Koelman's Theorem, these generate $I_\Delta$ as $\Delta$ has
 6 vertices.
 \hfill \raisebox{-3pt}{\epsfxsize=11pt\epsffile{figures/QED.eps}}
\end{example}

\section{Linear Projections}
Linear projections are key to the relationship between the toric varieties
$Y_\Delta$ introduced in Section 1 and Krasauskas's toric 
patches\index{toric patches}~\cite{Kr02}.
In geometric modeling, linear projections are encoded in the language
of weights and control points.
In algebraic geometry, linear projections provide the link between 
projective toric varieties $Y_{\mathcal{A}}$ and general rational 
varieties\index{rational variety}. 

Given $1+\ell$ vectors
${\bf p}_0, {\bf p}_1, \dotsc, {\bf p}_\ell \in\mathbb{C}^{1+k}$, we have the
linear map $\mathbb{C}^{1+\ell}\to\mathbb{C}^{1+k}$
 \begin{equation}\label{E:control}
  {\bf x}\ =\ (x_0,x_1,\dotsc,x_\ell)\ \longmapsto\ 
   x_0{\bf p}_0 + x_1{\bf p}_1 + \dotsb + x_\ell{\bf p}_\ell \ \in\ 
   \mathbb{C}^{1+k}\,,
 \end{equation}
represented by the matrix whose columns are the vectors ${\bf p}_i$.
Let $E:=\{{\bf x}\in\mathbb{C}^{1+\ell}\mid 0=\sum_i x_i{\bf p}_i\}$
be the kernel of this map.

Let $\mathbb{P}(E)$ be the linear subspace of $\mathbb{P}^\ell$ corresponding
to $E$.
Then~\eqref{E:control} induces a map $\pi$ from the difference 
$\mathbb{P}^\ell-\mathbb{P}(E)$ to $\mathbb{P}^k$, called a {\it linear
projection} with center of projection $\mathbb{P}(E)$
(or central projection from $\mathbb{P}(E)$).
We write
\[
   \pi\ \colon\ \mathbb{P}^\ell\ {-} {-} {-}{\to}\ \mathbb{P}^k\,.
\]
(In algebraic geometry, a broken arrow is used to represent such a 
{\it rational map}---a function that is not defined on all of $\mathbb{P}^\ell$.)
The {\it control points}\index{control point} ${\bf b}_0,\dotsc,{\bf b}_\ell$
of this projection are the images in $\mathbb{P}^k$ of the vectors ${\bf p}_i$.

Given a subvariety $Y\subset\mathbb{P}^\ell$ that does not meet the center
$\mathbb{P}(E)$, the linear projection $\pi$ restricts to give a map 
$\pi\colon Y\to\mathbb{P}^k$.
Points where $\mathbb{P}(E)$ meets $Y$ are called 
{\it basepoints}\index{basepoints} of the
projection $Y\,{-}{\to}\;\mathbb{P}^k$.

\begin{example}\label{Ex:toric_proj}
 A projective toric variety $Y_\mathcal{A}\subset\mathbb{P}^k$ as defined in
 Section~1 is the image of the projective toric variety
 $Y_\Delta\subset\mathbb{P}^\ell$, where $\Delta$ is the convex hull of the
 exponent vectors $\mathcal{A}\subset\mathbb{Z}^n$.
 Recall that the coordinates of $\mathbb{P}^k$ are naturally indexed by the
 elements of $\mathcal{A}$ and those of $\mathbb{P}^\ell$ by
 $\Delta\cap\mathbb{Z}^n$. 
 The projection simply `forgets' the coordinates of points whose index does not
 lie in $\mathcal{A}$.
 What is not immediate from the definitions is that the projection
 \begin{equation}\label{E:toric_proj}
   Y_\Delta\ \longrightarrow\ Y_{\mathcal{A}}
 \end{equation}
 has no basepoints.

 To see this, note that the center of this projection is defined by the vanishing
 of all coordinates of $\mathbb{P}^\ell$ indexed by elements of $\mathcal{A}$. 
 Any point ${\bf m}\in\Delta$ that is not a vertex is a positive rational
 combination of the vertices ${\bf v}$ of $\Delta$.
 That is, there are positive integers $d_{\bf m}$ and $d_{\bf v}$ for each
 vertex ${\bf v}$ of $\Delta$ such that 
\[
   d_{\bf m}\cdot {\bf m}\ =\ \sum_{\bf v} d_{\bf v}\cdot {\bf v}\,,
\]
 and so we have the binomial in the toric ideal $I_\Delta$
\[
   x_{\bf m}^{d_{\bf m}}\ -\ \prod_{\bf v} x_{\bf v}^{d_{\bf v}}\ 
\]

 In particular, if a point ${\bf x}\in Y_\Delta$ has a nonvanishing 
 ${\bf m}$th coordinate ($x_{\bf m}\neq 0$), then some vertex coordinates 
 $x_{\bf v}$ must also be nonvanishing.
 This shows that the map~\eqref{E:toric_proj} has no basepoints.

 This discussion also shows that the only basis vectors 
 $[0,\dotsc,0,1,0,\dotsc,0]$ contained in a projective toric variety
 $Y_{\mathcal{A}}$ are those indexed by the extreme points of
 $\mathcal{A}$---the vertices of the convex hull of $\mathcal{A}$.
 \hfill \raisebox{-3pt}{\epsfxsize=11pt\epsffile{figures/QED.eps}}
\end{example}

\begin{example}\label{E:ratcurve}
 Let ${\bf p}_0, {\bf p}_1, \dotsc, {\bf p}_n$ be vectors in $\mathbb{C}^{1+k}$. 
 Then the image $Z$ of the rational normal curve $Y_{[n]}$ of
 Example~\ref{E:3examples} under the corresponding linear projection is
 parametrized by
 \begin{equation}\label{E:ratCurve}
   \mathbb{P}^1\ \ni\ [s,t]\ \longmapsto\ 
    s^n{\bf p}_0 + s^{n-1}t{\bf p}_1+ \dotsb + t^n {\bf p}_n\ 
    \in\ \mathbb{P}^k\,.
 \end{equation}
 The map $Y_{[n]}\to Z$ has a basepoint at $[s,t]\in\mathbb{P}^1$
 when the sum in~\eqref{E:ratCurve} vanishes.
 Since each component of the sum is a homogeneous polynomial in $s,t$ of degree
 $n$, this implies that these $1+k$ polynomials share a common factor.
 When the polynomials have no common factor, $Z$ is a rational curve of degree
 $n$. 

 We consider an example of this when $n=3$ and $k=2$.
 Let $(1,-1,-1)$, $(1,-3,-1)$, $(1,-1,3)$, and $(1,1,-1)\in\mathbb{C}^3$ be the
 vectors ${\bf p}_0,\dotsc,{\bf p}_3$ which determine a linear projection
 $\mathbb{P}^3\,{-}{\to}\,\mathbb{P}^2$. 
 Then the image curve $Z$ of the toric variety $Y_{[3]}$ is given parametrically
 as  
\begin{eqnarray*}
   z_0&=& s^3+s^2t+st^2+t^3, \\
   z_1&=& -s^3-3s^2t-st^2+t^3, \ \ \textrm{and}\\
   z_2&=&  -s^3-s^2t+3st^2-t^3\,.
\end{eqnarray*}
 If we set $x= z_1/z_0$ and $y= z_2/z_0$ to be coordinates for the principal
 affine part of $\mathbb{P}^2$, then this has 
 implicit equation 
$$
  y^2(x-1)+2yx+x^2+x^3\ =\ 0\,.
$$
 We plot the control points\index{control point} ${\bf b}_i$ and the curve in
 Figure~\ref{F:cubic}. 
\begin{figure}[htb]
\[
  \begin{picture}(150,120)
   \put(0,0){\epsfysize=120pt\epsffile{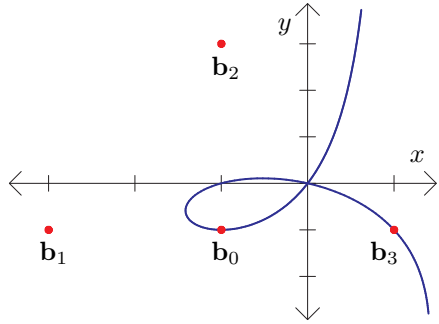}}
   \put(  0,23){${\bf b}_1$}
   \put( 65,23){${\bf b}_0$}
   \put(125,23){${\bf b}_3$}
   \put( 65,93){${\bf b}_2$}
   \put( 140, 60){$x$}
   \put(  90,110){$y$}
  \end{picture}
\]
\caption{A cubic curve}\label{F:cubic}
\end{figure}
\end{example}

\section{Rational Varieties}

Example~\ref{Ex:toric_proj} shows how the toric varieties $Y_{\mathcal{A}}$ and
$Y_\Delta$ are related via special linear projections
and Example~\ref{E:ratcurve} shows how rational curves are related to the
rational normal curve.
More general linear projections give 
rational varieties\index{rational variety}, which are varieties parametrized by
some collection of polynomials.

\begin{definition}
 A {\it rational variety} $Z\subset\mathbb{P}^k$ is the image of a projective
 toric variety\index{toric variety} $Y_{\mathcal{A}}$ under a linear
 projection. 
 The composition
\[
   (\mathbb{C}^*)^n \ \ \ 
    {\epsffile{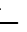}}%
    {\xrightarrow{\varphi_{\mathcal{A}}\ }}\ \ 
    Y_{\mathcal{A}}\ \ {-} {-} {-}{\to}\ \ Z
\]
 endows a rational variety $Z$ with a rational parametrization by polynomials
 whose monomials have exponent vectors in $\mathcal{A}$  and this
 parametrization is defined for (almost all) points in
 $(\mathbb{C}^*)^n$. 
\end{definition}

\begin{remark}
 The class of rational varieties is strictly larger than that of toric
 varieties.
 For example, the quartic rational plane curve whose rational parametrization
 and picture is shown 
 below is not a toric variety---its 
 three singular points prevent it from containing a dense torus.\smallskip

\ \ 
\begin{minipage}{1.4in}
 \epsfxsize=1.2in\epsffile{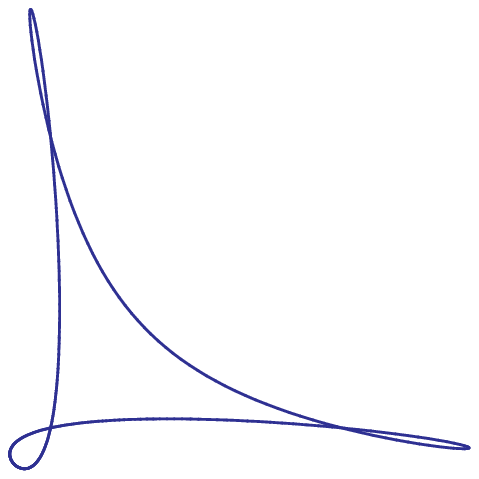}
 \end{minipage}
\begin{minipage}{2.7in}
   \begin{eqnarray*}
    x &=& t^4+7st^3+9s^2t^2-7s^3t-10s^4\\
    y &=& t^4-7st^3+9s^2t^2+7s^3t-10s^4\rule{0pt}{15pt}\\
    z &=& 3t^4-11s^2t^2+80s^4\rule{0pt}{15pt}\\
   \end{eqnarray*}
 \end{minipage}
\end{remark}\bigskip

This class of rational varieties contains 
the closures of the images of B\'ezier curves\index{B\'ezier curve}, triangular
B\'ezier patches\index{B\'ezier patch}, tensor product surfaces, and Krasauskas's toric
patches\index{toric patches}~\cite{Kr02}.
We give another example based upon the hexagon of Example~\ref{E:hex}.

\begin{example}\label{E:Hex_sur}
 Consider a projection $\mathbb{P}^6\,{-}{\to}\,\mathbb{P}^3$ where the 
 points ${\bf p}_i$ corresponding to the vertices  
  $\genfrac{[}{]}{0pt}{1}{1}{0},\genfrac{[}{]}{0pt}{1}{1}{1},
   \genfrac{[}{]}{0pt}{1}{0}{1},\genfrac{[}{]}{0pt}{1}{-1}{\ 0},
   \genfrac{[}{]}{0pt}{1}{-1}{-1},\genfrac{[}{]}{0pt}{1}{\ 0}{-1}$ 
 of the hexagon are 
 the following points in $\mathbb{C}^4$ taken in order:
\[
   (1,1,0,0), \ (1,1,1,0), \ (1,0,1,0), \ (1,0,1,1), \ (1,0,0,1), \ (1,1,0,1)
\]
 and suppose that the center of the hexagon corresponds to the point
 $(1,-1,-1,-1)$. 
 In  coordinates $[w,x,y,z]$ for $\mathbb{P}^3$, this has the rational
 parametrization.
 \begin{eqnarray*}
    w&=& 1+s+st+t+s^{-1}+s^{-1}t^{-1}+t^{-1}\\
    x&=&-1+s+st                      +t^{-1}\\ 
    y&=&-1  +st+t+s^{-1}                     \\
    z&=&-1       +s^{-1}+s^{-1}t^{-1}+t^{-1}
 \end{eqnarray*}
 Here are two views of (part of) the resulting rational surface and the axes.
\[
  \epsfxsize=1.9in\epsffile{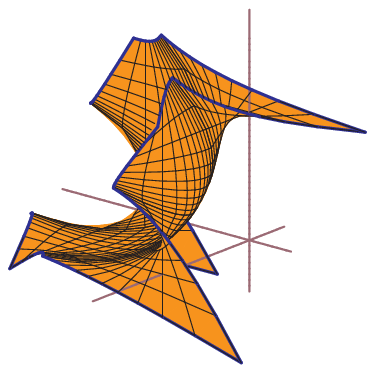}\qquad\qquad
  \epsfxsize=1.9in\epsffile{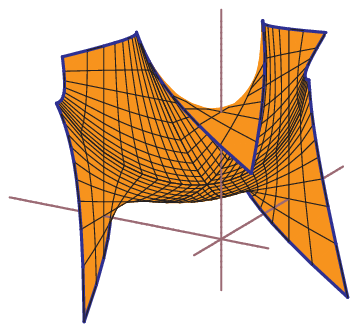}
\]

 As a subset of $\mathbb{P}^3$, this is defined by a the vanishing of a single
 homogeneous polynomial of degree 6 with 72 terms\medskip 

\noindent\hspace{2.1pt}
 $112w^6
  -240(x+y+z)w^5
  +296(xy+xz+yz)w^4+216(x^2+y^2+z^2)w^4$

\noindent\hspace{2.1pt}
 $-92(x^3+y^3+z^3)w^3
  -124(x^2y+xy^2+x^2z+xz^2+y^2z+yz^2)w^3
  -568xyzw^3$

\noindent\hspace{2.1pt}
 $+4(x^4+y^4+z^4)w^2
  +70(x^3y+xy^3+x^3z+xz^3+y^3z+yz^3)w^2$

\noindent\hspace{2.1pt}
 $-125(x^2y^2+x^2z^2+y^2z^2)w^2
  +272(x^2yz+xy^2z+xyz^2)w^2$

\noindent\hspace{2.1pt}
 $-2(x^4y+xy^4+x^4z+xz^4+y^4z+yz^4)w
  -141(x^3yz+xy^3z+xyz^3)w$

\noindent\hspace{2.1pt}
 $+35(x^3y^2+x^2y^3+x^3z^2+x^2z^3+y^3z^2+y^2z^3)w
  -7(x^2y^2z+x^2yz^2+xy^2z^2)w$

\noindent\hspace{2.1pt}
 $+5(x^4yz+xy^4z+xyz^4)
  +19(x^3y^2z+x^3yz^2+x^2y^3z+x^2yz^3+xy^3z^2+xy^2z^3)$

\noindent\hspace{2.1pt}
 $-50x^2y^2z^2
  -13(x^3y^3+x^3z^3+y^3z^3)
  -2(x^4y^2+x^2y^4+x^4z^2+x^2z^4+y^4z^2+y^2z^4)$\medskip 

\noindent
The symmetry of this polynomial in the variables $x,y,z$
is due to the symmetry of the hexagon and of the control 
points. 
 \hfill \raisebox{-3pt}{\epsfxsize=11pt\epsffile{figures/QED.eps}}
\end{example}

In these two examples, the toric ideals of the varieties $Y_{\mathcal{A}}$
were generated by quadratic binomials, while the resulting
rational varieties were hypersurfaces defined by polynomials of degrees 4 and 6
respectively. 
These examples show that the ideal of a rational variety may be rather
complicated.
Nevertheless, this ideal can be computed quite reasonably either from the
original toric ideal and the projection or from the resulting parametrization. 
(See Sections 3.2 and 3.3 of~\cite{CLO97} for details.)

\section{Implicit Degree of a Toric Variety}

The (implicit) degree of a hypersurface (e.g. planar curve or a surface in
$\mathbb{P}^3$) is the degree of its implicit equation.
Similarly, the degree of a rational curve 
${\bf f}\colon\mathbb{P}^1\to\mathbb{P}^k$
is the degree of the components of its parametrization ${\bf f}$.
Projective varieties with greater dimension and codimension also have a degree
that is well-behaved under linear projection, and this degree is readily
determined for toric varieties. 

\begin{definition}
 Let $X\subset\mathbb{P}^\ell$ be an algebraic variety of dimension $n$.
 The {\it degree} of $X$, $\textrm{deg}(X)$, is the number of points common to
 $X$ and to a general linear subspace $L$ of dimension $\ell-n$.
 Such a linear subspace is defined by $n$ linear equations and so
 the degree of $X$ is also the number of (complex) solutions to $n$
 general linear equations on $X$.
\end{definition}

\begin{remark}
 This notion of degree agrees with the usual notions for hypersurfaces
 and for rational curves.
 Suppose that $X\subset\mathbb{P}^\ell$ is a hypersurface with implicit equation
 $f=0$ where $f$ has degree $d$.
 Restricting $f$ to a line $L$ and identifying $L$ with
 $\mathbb{P}^1$ gives a polynomial of degree $d$ on $\mathbb{P}^1$ whose zeroes
 correspond to the points of $X\cap L$.
 If $L$ is general then there will be $d$ distinct roots of this polynomial,
 showing the equality $\deg(X)=\deg(f)$, as $\deg(X)$ is
 $\# L\cap X$ which equals $d= \deg(f)$.

 Similarly, suppose that $X\subset\mathbb{P}^\ell$ is a rational curve of degree
 $d$. 
 Then it has a parametrization 
 ${\bf f}\colon\mathbb{P}^1\to X\subset\mathbb{P}^\ell$
 where the components of ${\bf f}$ are homogeneous polynomials of degree $d$.
 A hyperplane $L$ in $\mathbb{P}^\ell$ is defined by a single linear equation
 $\Lambda({\bf x})=0$.
 Then the points of $X\cap L$ correspond to the zeroes of the polynomial
 $\Lambda({\bf f})$, which has degree $d$ (when $L$ does not contain the image
 of ${\bf f}$).
 \hfill \raisebox{-3pt}{\epsfxsize=11pt\epsffile{figures/QED.eps}}
\end{remark}

If $\pi\colon\mathbb{P}^\ell\,{-}{\to}\,\mathbb{P}^k$ is a surjective linear
projection with center $\mathbb{P}(E)$, then the inverse image of a linear
subspace $K\subset\mathbb{P}^k$ of dimension $k-n$ is a linear subspace $L$ of 
$\mathbb{P}^\ell$ of dimension $\ell-n$ that contains the center
$\mathbb{P}(E)$.  
This implies that the degree of a projective variety is reasonably well-behaved
under linear projection.
We give the precise statement.

\begin{theorem}
  Let $\pi\colon\mathbb{P}^\ell\,{-}{\to}\,\mathbb{P}^k$ be a linear projection  
  with center $\mathbb{P}(E)$ and $Y\subset\mathbb{P}^\ell$.
  Suppose that $\pi(Y)$ has the same dimension as does $Y$.
  Then
\[
   \deg(\pi(Y))\  \leq\ \deg(Y)\,,
\]
 with equality when the map $\pi\colon Y\to\pi(Y)$ has no basepoints
 and is one-to-one.
 
 These conditions are satisfied for a general linear projection if $\dim Y<k$.
\end{theorem}

Thus rational varieties $Z$ that have an injective (one-to-one) parametrization
given by a map $\pi\colon Y_{\mathcal{A}}\to Z$ with no basepoints will have the
same degree as the projective toric variety $Y_{\mathcal{A}}$.
This degree is nicely expressed in terms of  the convex hull $\Delta$ of the
exponent vectors $\mathcal{A}$.

\begin{theorem}
 The implicit degree of a toric variety $Y_{\mathcal{A}}$ is
\[
     n!\mbox{Vol}(\Delta)\,,
\]
  where $\mbox{Vol}(\Delta)$ is the usual Euclidean volume of the
  $n$-dimensional polytope $\Delta$.
\end{theorem}

Thus the degree of a rational variety $Z$ parametrized by polynomials whose
monomials have exponents from a set $\mathcal{A}$ whose convex hull is $\Delta$
is at most $n!\mbox{Vol}(\Delta)$, with equality when the parametrization
$\pi\colon Y_{\mathcal{A}}\to Z$ has no basepoints and is one-to-one
(injective). 
This is what we saw in the examples of Section~4.
The polytope of the quartic curve is a line segment of length, and hence
volume, 4, while the hexagon whose corresponding monomials parameterize the
rational surface of Example~\ref{E:Hex_sur}
has area 3, and $2!\cdot 3=6$, which is the degree of its implicit
equation.\smallskip 

This determination of the degree of the toric variety $Y_{\mathcal{A}}$ is an
important result due to Kouchnirenko~\cite{BKK} concerning the solutions of
{\it sparse equations}.  
One of the equivalent definitions of the degree of $Y_{\mathcal{A}}$ 
is the number of solutions to $n$ ($=\dim(Y_{\mathcal{A}})$) equations on
$Y_{\mathcal{A}}\subset\mathbb{P}^\ell$.
Under the parametrization $\varphi_{\mathcal{A}}$~\eqref{E:calA} of
$Y_{\mathcal{A}}$,
these linear equations become {\it Laurent polynomials} on $(\mathbb{C}^*)^n$
whose monomials have exponent vectors in $\mathcal{A}$.
Thus the degree of $Y_{\mathcal{A}}$ is equal to the number of solutions
in $(\mathbb{C}^*)^n$ to $n$ general Laurent polynomials  whose monomials have
exponent vectors in $\mathcal{A}$.

This result of Kouchnirenko was generalized by Bernstein~\cite{Be75} who
determined the number of solutions in 
$(\mathbb{C}^*)^n$ to $n$ general Laurent polynomials with possibly different
sets of exponent vectors.
In that, the r\^{o}le of the volume is played by the mixed volume.
For more, see the contribution of Rojas~\cite{Ro03} to these proceedings.

\section{The Real Part of a Toric Variety}

B\'ezier curves and surface patches in geometric modeling are parametrizations
of some of the real part of a rational variety.\index{toric variety!real}
We discuss the real part of a toric variety and of rational varieties, with
respect to their usual real structure.
Some toric varieties admit exotic real structures, a topic covered in 
the article by Delaunay~\cite{De03} that also appears in this volume.

\begin{definition}
 The (standard) real part of a toric variety is defined by 
 replacing the complex numbers $\mathbb{C}$ by the real numbers $\mathbb{R}$
 everywhere in the given definitions.
\end{definition}

For example, consider the projective toric variety $Y_{\mathcal{A}}$, defined as
a subset of projective space $\mathbb{P}^\ell$ by the toric ideal
$I_{\mathcal{A}}$ (equivalently, as the closure of the image of
$\varphi_\mathcal{A}$~\eqref{E:calA}). 
Then the real part $Y_{\mathcal{A}}(\mathbb{R})$ of $Y_{\mathcal{A}}$ is the 
intersection of $Y_{\mathcal{A}}$ with $\mathbb{RP}^\ell$,
that is, the subset of $\mathbb{RP}^\ell$ defined by the toric ideal.
Recall that $I_{\mathcal{A}}$ is generated
by binomials ${\bf x}^{\bf u}-{\bf x}^{\bf v}$,  which are real polynomials.

Suppose that we have a linear projection 
$\pi\colon\mathbb{P}^\ell\,{-}{\to}\,\mathbb{P}^k$ defined by real
points ${\bf p}_i\in\mathbb{R}^{1+k}$.
Then the rational variety $Z$ (the image of
$Y_{\mathcal{A}}$ under $\pi$) has ideal $I(Z)$ generated by real polynomials.
The real part $Z(\mathbb{R})$ of $Z$ is the subset of $\mathbb{RP}^k$ defined by
the ideal $I(Z)$. 
This again coincides with the intersection of $Z$ with $\mathbb{RP}^k$.
All pictures in this tutorial arise in this fashion. 
When the map $\pi$ has no basepoints and 
$\pi\colon Y_{\mathcal{A}}\to Z$ is one-to-one at almost all points of $Z$, then
$\pi(Y_{\mathcal{A}}(\mathbb{R}))=Z(\mathbb{R})$.

The reason for this is that when $x\in\mathbb{RP}^k$, the points
in $\pi^{-1}(x)\cap Y_{\mathcal{A}}$ are the solution to a system of equations
with real coefficients.  
Since the map $\pi$ is one-to-one on $Y_{\mathcal{A}}$, this system has a single
solution that is necessarily real.
If $\pi$ is not one-to-one, then we may have 
$\pi(Y_{\mathcal{A}}(\mathbb{R}))\subsetneq Z(\mathbb{R})$.
For example, when $\mathcal{A}$ is a line segment of length 2,
$Y_{\mathcal{A}}$ is the parabola 
$\{[1,x,x^2]\mid x\in\mathbb{C}\}\subset\mathbb{P}^2$. 
The projection to $\mathbb{P}^1$ omitting the second coordinate (which is
basepoint-free) is the two-to-one map
\[
   \mathbb{C}\ni x\mapsto [1,x^2]\in\mathbb{P}^1
\]
whose restriction to $\mathbb{R}$ has image the nonnegative part of the
real toric variety $\mathbb{RP}^1$.

This description does little to aid our intuition about the real part of a toric
variety or a rational variety.
We obtain a more concrete picture of the real points of a toric variety
$Y_{\mathcal{A}}$ 
and an alternative construction of $Y_{\mathcal{A}}(\mathbb{R})$ 
if we first describe the real
points of an abstract toric variety $X_\Sigma$.
For this, we recall the definition of the abstract toric variety $X_\Sigma$
associated to a fan $\Sigma$, as described by Cox~\cite{Co02}.

Let $\sigma\subset\mathbb{R}^n$ be a strongly convex
rational polyhedral cone with dual cone $\sigma^\vee\subset\mathbb{R}^n$.
Lattice points ${\bf m}\in\sigma^\vee\cap\mathbb{Z}^n$ are exponent vectors of
Laurent monomials ${\bf t}^{\bf m}$ defined on $(\mathbb{C}^*)^n$.
The affine toric variety $U_\sigma$ corresponding to $\sigma$ is constructed by
first choosing a finite generating set  
${\bf m}_1,{\bf m}_2,\dotsc,{\bf m}_\ell$ of the additive semigroup 
$\sigma^\vee\cap\mathbb{Z}^n$.
These define the map
 \begin{eqnarray*}
  \varphi\ \colon\ (\mathbb{C}^*)^n&\longrightarrow&\mathbb{C}^\ell\\
   {\bf t}&\longmapsto&
  ({\bf t}^{{\bf m}_1},{\bf t}^{{\bf m}_2},\dotsc,{\bf t}^{{\bf m}_\ell})\ ,
 \end{eqnarray*}
and we set $U_\sigma$ to be the closure of the image of this map.
The real part $U_\sigma(\mathbb{R})$ of this affine toric variety is simply the
intersection of $U_\sigma$ with $\mathbb{R}^\ell$.

The intersection $\sigma\cap\tau$ of two cones $\sigma,\tau$ in a fan $\Sigma$
is a face of each cone and $U_{\sigma\cap\tau}$ is naturally a subset of
both $U_\sigma$ and $U_\tau$.
The toric variety $X_\Sigma$ is obtained by gluing together $U_\sigma$ and
$U_\tau$ along their common subset $U_{\sigma\cap\tau}$, for all cones 
$\sigma,\tau$ in $\Sigma$.
The real part $X_\Sigma(\mathbb{R})$ of $X_\Sigma$ is similarly obtained by
piecing together the real parts $U_\sigma(\mathbb{R})$ and $U_\tau(\mathbb{R})$
along their common subset $U_{\sigma\cap\tau}(\mathbb{R})$, for all
$\sigma,\tau$ in $\Sigma$.

Since the origin ${\bf 0}\in\mathbb{R}^n$ lies in $\Sigma$ and 
${\bf 0}^\vee=\mathbb{R}^n$, the affine toric variety $U_{{\bf 0}^\vee}$ is the 
torus $(\mathbb{C}^*)^n$, which is a common subset of each $U_\sigma$. 
This torus is dense in the toric variety $X_\Sigma$ and it acts on $X_\Sigma$.
Similarly, the torus $(\mathbb{R}^*)^n$ is dense in $X_\Sigma(\mathbb{R})$ and
it acts on $X_\Sigma(\mathbb{R})$. 
This torus $(\mathbb{R}^*)^n$ has $2^n$ components called {\it orthants}, 
each identified by the sign vector $\varepsilon\in\{\pm1\}^n$ recording the
signs of coordinates of points in that component.
The identity component is the orthant containing the identity, and it has sign
vector $(1,1,\dotsc,1)$.
Write $\mathbb{R}_>^n$ for this identity component.

\begin{definition}
 The {\it non-negative part} $X_\geq$ of a toric variety $X$ is the closure 
 (in $X(\mathbb{R})$) of the identity component $\mathbb{R}_>^n$ of
 $(\mathbb{R}^*)^n$. 
 The {\it boundary} of $X_\geq$ is defined to be the difference 
 $X_\geq-\mathbb{R}_>^n$.
\end{definition}

We could also consider the closures of other components of the torus
$(\mathbb{R}^*)^n$, obtaining $2^n$ other pieces analogous to this non-negative
part $X_\geq$. 
Since the component of $(\mathbb{R}^*)^n$ having sign vector $\varepsilon$ is
simply $\varepsilon\cdot\mathbb{R}_>^n$, these other pieces are
translates of $X_\geq$ by the appropriate sign vector, and hence all are
isomorphic. 
Since $X(\mathbb{R})$ is the closure of $(\mathbb{R}^*)^n$ and each 
piece $\varepsilon\cdot X_\geq$ is the closure of the orthant
$\varepsilon\cdot\mathbb{R}_>^n$, we obtain a concrete picture of
$X(\mathbb{R})$: it is pieced
together from $2^n$ copies of this non-negative part $X_\geq$ glued together
along their common boundaries.

The non-negative part of the toric variety $Y_\mathcal{A}$ is simply the
intersection of $Y_\mathcal{A}$ with the non-negative part of the ambient
projective space $\mathbb{P}^\ell$, those points with non-negative homogeneous
coordinates 
\[
   \{ [x_0,x_1,\dotsc,x_\ell]\mid x_i\geq 0\}\,.
\]
The boundary of $(Y_\mathcal{A})_\geq$ is its intersection with the coordinate
hyperplanes, which are defined by the vanishing of at least one homogeneous
coordinate.

\begin{example}
 The surface of Example~\ref{E:Hex_sur} is the image of the toric variety
 $Y_\Delta$, where $\Delta$ is the hexagon of Example~\ref{E:hex}.
 Figure~\ref{F:hex} shows the image of the non-negative part of $Y_\Delta$.
 The control points are the spheres (dots) and the boundary consists of the thickened
 lines.
\begin{figure}[htb]
\[
   \epsfxsize=2.7in\epsffile{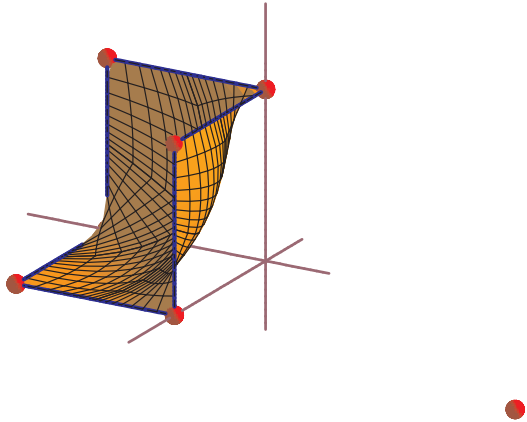}
\]
 \caption{Hexagonal toric patch}\label{F:hex}
\end{figure}
The six control points corresponding to the vertices of the hexagon lie on the
non-negative part of $Y_\Delta$. 
The seventh control point corresponding to the center of the hexagon appears in the
lower right.
It lies in the octant opposite to the non-negative part of $Y_\Delta$ and 
causes $Y_\Delta$ to `bulge' towards the origin.
\end{example}

\section{The Double Pillow}

We devote this section to the construction of the toric variety
$X_\Sigma(\mathbb{R})$ for a single example, where $\Sigma$ is the normal
fan\index{fan!normal, of a polytope} of 
the cross polytope $\lozenge\subset\mathbb{R}^2$.
As remarked in Section 13 of~\cite{Co02}, $X_\Sigma\simeq Y_\lozenge$ as
$\lozenge$ is 2-dimensional.
Krasauskas~\cite{Kr02} calls the corresponding toric surface the `pillow with
antennae'. 
We display $\lozenge$ together with its normal fan $\Sigma$, with one of its
full-dimensional cones shaded.
 \begin{equation}\label{E:loznege}
   \raisebox{-0.68in}{
      \begin{picture}(105,112)
       \put(0,0){\epsfysize=1.44in\epsffile{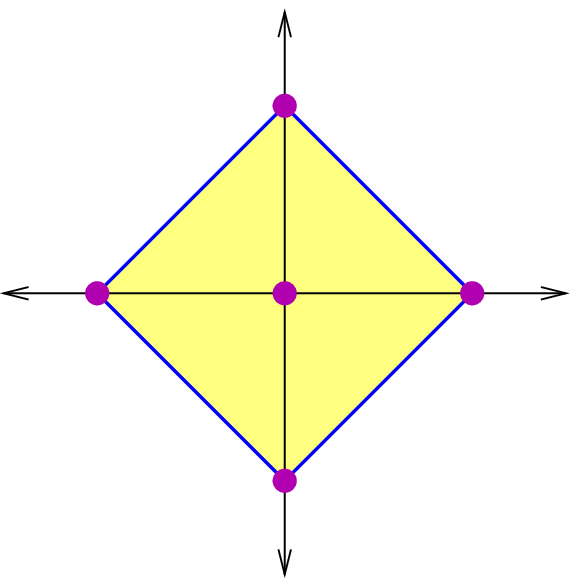}}
      \end{picture}}
      \qquad\qquad
   \raisebox{-0.68in}{
      \begin{picture}(125,112)
       \put(0,0){\epsfysize=1.48in\epsffile{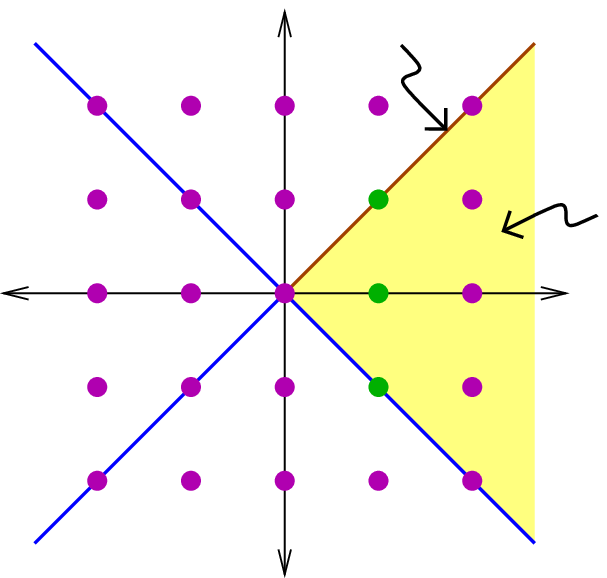}}
       \put(69,101){$\tau$}
       \put(115,66){$\sigma$}
      \end{picture}}
 \end{equation}

Each full-dimensional cone $\sigma$ is self-dual and they are all isomorphic.
Thus $Y_\lozenge(\mathbb{R})$ ($=X_\Sigma(\mathbb{R})$) is obtained by gluing
together four isomorphic affine toric varieties $U_\sigma(\mathbb{R})$, as
$\sigma$ ranges over the 2-dimensional cones in $\Sigma$.
A complete
picture of the gluing involves the affine varieties $U_\tau(\mathbb{R})$, where
$\tau$ is one of the rays of $\Sigma$.
We next describe these two toric varieties $U_\sigma(\mathbb{R})$ and
$U_\tau(\mathbb{R})$, for $\sigma$ a 2-dimensional cone and $\tau$ a ray of
$\Sigma$. 

Let $\sigma$ be the shaded cone in~\eqref{E:loznege}.
Since $\sigma=\sigma^\vee$, we see that $\sigma^\vee\cap\mathbb{Z}^2$ is
minimally generated by the vectors $(1,-1), (1,0)$, and $(1,1)$, and so
$U_\sigma(\mathbb{R})$ is the closure in $\mathbb{R}^3$ of the image of the map
\[
  \varphi\ \colon\ (s,t)\ \longmapsto\ (st^{-1}, st, s)\,,
\]
which is defined by the equation $xy=z^2$ (where $(x,y,z)$ are the coordinates
for $\mathbb{R}^3$).
This is a right circular cone in $\mathbb{R}^3$, which we display 
below at left.
\[
  \begin{picture}(144,115)(0,-15)
  \put(0,0){\epsfxsize=2in\epsffile{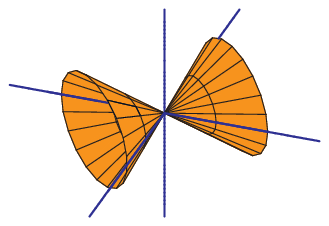}}
   \put(108,88){$y$}
   \put( 65,87){$z$}
   \put(135,41){$x$}
   \put(62,-15){$U_\sigma(\mathbb{R})$}
  \end{picture}
   \qquad\qquad
  \begin{picture}(144,115)(0,-15)
   \put(0,0){\epsfxsize=2in\epsffile{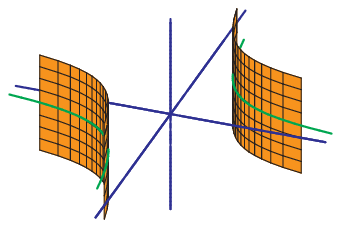}}
   \put(110,90){$y$}
   \put( 65,87){$z$}
   \put(140,28){$x$}
   \put(64,-15){$U_\tau(\mathbb{R})$}
  \end{picture}
\]

Let $\tau$ be the ray generated by $(1,1)$, which is a face of $\sigma$.
Then $\tau^\vee$ is the half-space $\{(u,v)\in\mathbb{R}^2\mid u+v\geq 0\}$, 
which is the union of both  2-dimensional cones in $\Sigma$ containing $\tau$.
Since $\tau^\vee\cap\mathbb{Z}^2$ has generators $(1,-1), (-1,1)$, and $(1,0)$,
we see that $U_\tau(\mathbb{R})$ is the closure in $\mathbb{R}^3$ of the image
of the map
\[
  \varphi\ \colon\ (s,t)\ \longmapsto\ (st^{-1}, s^{-1}t, s)\,,
\]
which has equation $xy=1$.
This is the cylinder with base the hyperbola $xy=1$, which is shown 
above at right.

We describe the gluing.
We know that $U_\tau(\mathbb{R})\subset U_\sigma(\mathbb{R})$ and they both
contain the torus $(\mathbb{R}^*)^2$.
This common torus is their intersection with the complement of the
coordinate planes, $xyz\neq 0$, and their boundaries are their intersections 
with the coordinate planes.
The boundary of the cylinder is the curve $z=0$ and $xy=1$, which is defined by
$s=0$ and displayed on the picture of $U_\tau(\mathbb{R})$.
Also, $t\neq 0$ on the cylinder.
The boundary of the cone is the union of the $x$ and $y$ axes.  
Since $t^2=y/x$ on the cone, the locus where $t=0$ is the $x$ axis.
Thus $U_\tau(\mathbb{R})$ is naturally identified with the complement of the $x$ 
axis in $U_\sigma(\mathbb{R})$ where the curve  $z=0, xy=1$ in
$U_\tau(\mathbb{R})$ is identified with the $y$-axis in $U_\sigma(\mathbb{R})$.

\begin{figure}[htb]
\[
  \epsfxsize=2.4in\epsffile{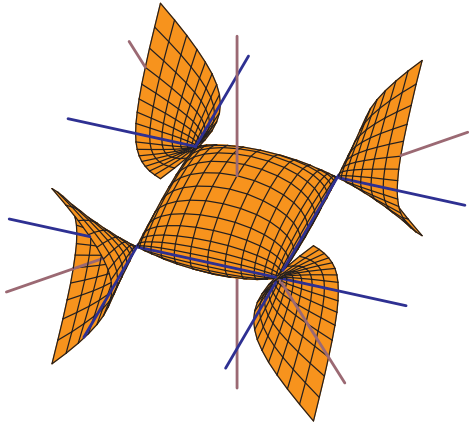}
\]
\caption{The Double Pillow.\label{F:pillow}}
\end{figure}

If $\tau'$ is the other ray defining $\sigma$, then
$U_{\tau'}(\mathbb{R})$ ($\simeq U_\tau(\mathbb{R})$) 
is identified with the complement of the $y$ axis in
$U_\sigma(\mathbb{R})$. 
A convincing understanding of this gluing procedure is obtained by considering
the rational surface $Z$ in $\mathbb{RP}^3$ which is the image of the toric
variety $Y_\lozenge(\mathbb{RP}^3)$ under the projection map given by the points
$(1,\pm1,0,0)$ and $(1,0,\pm1,0)$ associated to the vertices 
$(\pm1,0)$ and $(0,\pm1)$ of $\lozenge$,
and $(0,0,0,1)$ associated to its center.
We display this surface in Figure~\ref{F:pillow}.

This surface has the implicit equation
\[
  (x^2-y^2)^2 - 2x^2w^2 - 2y^2w^2 -16z^2w^2 + w^4\ =\ 0\,.
\]
and its dense torus has parametrization
\[
   [w,x,y,z]\ =\ 
     [{\textstyle s+t+\frac{1}{s}+\frac{1}{t},\  s-\frac{1}{s},\ 
      t-\frac{1}{t},\ 1}]\,.
\]
It has curves of self-intersection along the lines $x=\pm y$ in the plane
at infinity ($w=0$).
As the self-intersection is at infinity, this affine surface is a good
illustration of the toric variety $Y_\lozenge(\mathbb{R})$, and 
so we refer to this picture to describe $Y_\lozenge(\mathbb{R})$. 

This surface contains 4 lines $x\pm y=\pm1$ and 
their complement is the dense torus in $Y_\lozenge(\mathbb{R})$.
The complement of any three lines is the piece $U_\tau(\mathbb{R})$
corresponding to a ray $\tau$.
Each of the four singular points is a singular point of one cone
$U_\sigma(\mathbb{R})$, which is obtained by removing the two lines not meeting
that singular point.
Finally, the action of the group $\{(\pm1,\pm1)\}$ on
$Y_\lozenge(\mathbb{R})$ may also be seen from this picture.
Each singular point is fixed by this group.
The element $(-1,-1)$ sends $z\mapsto -z$, interchanging the top and bottom
halves of each piece, while the elements $(1,-1)$ and $(-1,1)$ interchange the
central `pillow' with the rest of $Y_\lozenge(\mathbb{R})$.
In this way, we see that $Y_\lozenge(\mathbb{R})$ is a `double pillow'.\smallskip

The non-negative part of $Y_\lozenge(\mathbb{R})$ is also readily determined
from this picture. 
The upper part of the middle pillow is the part of $Y_\lozenge(\mathbb{R})$
parametrized by $\mathbb{R}_>^2$, and so its closure is just a square, but with
singular corners obtained by cutting a cone into two pieces along a plane of
symmetry. 
In fact, the orthogonal projection to the $xy$ plane identifies this
non-negative part with the cross polytope $\lozenge$.
From the symmetry of this surface, we see that $Y_\lozenge(\mathbb{R})$ is
obtained by gluing four copies of cross polytope $\lozenge$ together along their
edges to form two pillows attached at their corners.
(The four `antennae' are actually the truncated corners of the second
pillow---projective geometry can play tricks on our affine intuition.)

\section{Linear Precision and the Algebraic Moment Map}

We observed  that the non-negative part of the toric
variety $Y_\lozenge$ can be identified with $\lozenge$. 
The non-negative part of any projective toric variety $Y_\mathcal{A}$ admits
an identification with the convex hull $\Delta $ of $\mathcal{A}$.
One way to realize this identification is through the moment 
map\index{moment map} and algebraic moment map of a toric variety $Y_\mathcal{A}\to\Delta$.

\begin{definition}
 Let $Y_\mathcal{A}\subset\mathbb{P}^\ell$ be a projective toric variety given
 by a collection of exponent vectors $\mathcal{A}\subset\mathbb{R}^n$
 with convex hull $\Delta$.
 The torus $(\mathbb{C}^*)^n$ acts on $\mathbb{P}^\ell$ and on $Y_\mathcal{A}$ 
 via the map $\varphi_\mathcal{A}$. 
 To such an action, symplectic geometry associates a {\it moment map}
 $\mu_\mathcal{A}\colon \mathcal{A}\to \mathbb{R}^n$, 
 \begin{equation}\label{E:moment}
  \begin{array}{rclcl}
    \mu_\mathcal{A}&\colon&Y_\mathcal{A}&\longrightarrow&\mathbb{R}^n\\
             &&{\bf x}&\longmapsto&{\displaystyle
              \frac{1}{\sum_{{\bf m}\in\mathcal{A}}|x_{\bf m}({\bf x})|^2}
             \sum_{{\bf m}\in\mathcal{A}}|x_{\bf m}({\bf x})|^2{\bf m}}\ .
  \end{array}
 \end{equation}
 While the restriction of coordinate function
 $x_{\bf m}$ on $\mathbb{P}^\ell$ to $Y_\mathcal{A}$ is not
 a well-defined function, the collection of these
 coordinate functions is well-defined up to a common scalar factor.
 It is a basic theorem of symplectic geometry that the image of the moment map 
 is the polytope $\Delta$ and the restriction of $\mu_\mathcal{A}$ to the non-negative
 part of $(Y_\mathcal{A})_\geq$ is a homeomorphism.
\end{definition}

 More useful to us is the following variant of $\mu_\mathcal{A}$, where we do not 
 square the coordinate functions,
 \begin{equation}\label{E:Amoment}
  \begin{array}{rclcl}
    \alpha_\mathcal{A}&\colon&Y_\mathcal{A}&\longrightarrow&\mathbb{R}^n\\
             &&{\bf x}&\longmapsto&{\displaystyle
              \frac{1}{\sum_{{\bf m}\in\mathcal{A}}|x_{\bf m}({\bf x})|}
             \sum_{{\bf m}\in\mathcal{A}}|x_{\bf m}({\bf x})|{\bf m}}\ .
  \end{array}
 \end{equation}
 This map is very similar to the moment map, and thus is often confused with the
 moment map\newfootnote{This confusion occurred in the published version of this
 manuscript.}

\begin{remark}\label{R:moment_proj}
 Suppose that 
 $\mathcal{A}=\{{\bf m}_0,{\bf m}_1,\dotsc,{\bf m}_\ell\}\subset\mathbb{R}^n$.
 We claim that on $(Y_\mathcal{A})_\geq$, the map~\eqref{E:Amoment}
 coincides with the linear projection 
 $\pi_\mathcal{A}\colon\mathbb{P}^\ell -\to \mathbb{P}^n$ defined by the points
 in the lift $\mathcal{A}^+$ of $\mathcal{A}$: 
\[
   (1,{\bf m}_0), (1,{\bf m}_1), \dotsc, (1,{\bf m}_\ell)\ \in\
   \mathbb{R}^{1+n}\ .
\]
 Indeed, we have
\[
  \pi_\mathcal{A}({\bf x})\ =\ 
     \pi_\mathcal{A}([x_0,x_1,\dotsc,x_\ell])\ =\ 
       \sum_{i=0}^{\ell} x_i [1,{\bf m}_i]\ =\ 
      \bigl[{\textstyle \sum_i x_i,\ \sum_i x_i{\bf m}_i}\bigr]\ .
\]
 If ${\bf x}$ lies in the non-negative part of the projective 
 toric variety $Y_{\mathcal{A}}$, then each coordinate $x_i$ of 
 ${\bf x}$ is non-negative with $x_i=x_{{\bf m}_i}$. 
 Since $\sum_ix_i>0$, this shows that 
 $\pi_\mathcal{A}({\bf x})=[1,\alpha_{\mathcal{A}}({\bf x})]$, and thus 
 the map~\eqref{E:Amoment} coincides with the projection $\pi_\mathcal{A}$
 on the non-negative part of $Y_{\mathcal{A}}$.

 It is for these reasons that we call the linear projection $\pi_\mathcal{A}$ the 
 {\it algebraic moment map}.
 \hfill \raisebox{-3pt}{\epsfxsize=11pt\epsffile{figures/QED.eps}}
\end{remark}

 This algebraic moment map shares an important property of the moment map. 

\begin{theorem}\label{T:moment}
  The non-negative part $(Y_\mathcal{A})_\geq$ of the toric variety
  $Y_\mathcal{A}$ is homeomorphic to the convex hull $\Delta$ of $\mathcal{A}$
  under the algebraic moment map. 
\end{theorem}

The nature of this homeomorphism is subtle.
If the polytope $\Delta$ is smooth (that is, the shortest integer vectors normal
to the faces that meet at a vertex always form a basis for $\mathbb{Z}^n$), then 
every point of $\Delta$ has a neighborhood in $\Delta$ homeomorphic to 
$\mathbb{R}^k\times \mathbb{R}^{n-k}_\geq$, and so we call $\Delta$ a 
{\it manifold with corners}.
In general, a polytope $\Delta$ is a manifold with `singular corners'.
It is this structure that is preserved by the homeomorphism of
Theorem~\ref{T:moment}. 
(For more on the algebraic moment map and the structure of $X_\geq$ as a manifold with
singular corners, see Section 4 of Fulton's book on toric
varieties~\cite{Fu93}, where he call $\alpha_\mathcal{A}$ the moment map.)\medskip 

Theorem~\ref{T:moment} explains why toric patches\index{toric patches} are of
interest in geometric modeling. 
Since the non-negative part of a projective toric surface is homeomorphic to a
polygon, any rational surface parametrized by that toric surface has a
non-negative part that is the image of that polygon.
In this way, we can obtain multi-sided surface patches from toric surfaces.
This theorem not only explains the geometry of such toric 
patches, but we use it to gain insight into parametrizations of toric
patches by the corresponding polytopes.

Let $\Delta$ be the convex hull of a set of exponent vectors $\mathcal{A}$.
By Theorem~\ref{T:moment}, $(Y_\mathcal{A})_\geq$ is homeomorphic to
$\Delta$, and so there exists a parametrization of $(Y_\mathcal{A})_\geq$ by 
$\Delta$ preserving their structures as manifolds with `singular corners'.
From the point of view of algebraic geometry, the most natural such
parametrization is the inverse of the algebraic moment map 
$\alpha_\mathcal{A}^{-1}\colon \Delta \to (Y_\mathcal{A})_\geq$.
This is also the most natural from the point of view of geometric modeling.

\begin{theorem}\label{T:precision}
 The coordinate functions of $\alpha_\mathcal{A}^{-1}$ have linear
 precision\index{linear precision}.
\end{theorem}

A collection of liinearly independent functions 
$\{f_{\bf m}\mid {\bf m}\in\mathcal{A}\}$  defined on the convex hull
$\Delta\subset\mathbb{R}^n$ of $\mathcal{A}$ has 
{\it linear precision} if, for any affine function $\Lambda$ defined on $\mathbb{R}^n$, 
 \begin{equation}\label{E:precision}
   \Lambda({\bf u})\ =\ 
    \sum_{{\bf m}\in\mathcal{A}} \Lambda({\bf m}) f_{\bf m}({\bf u})\
   \qquad\textrm{for all }{\bf u}\in \Delta\,.
 \end{equation}
Theorem~\ref{T:precision} follows immediately from this definition.
The functions $\{f_{\bf m}\mid {\bf m}\in\mathcal{A}\}$ define a map ${\bf f}$
from $\Delta$ to $\mathbb{P}^\ell$ in the natural coordinates of
$\mathbb{P}^\ell$ indexed by the exponent vectors in $\mathcal{A}$.
Then the right hand side of~\eqref{E:precision} is the result of applying the
linear function $\Lambda$ to the composition 
\[
   \Delta\ \xrightarrow{\ {\bf f}\ }\ (Y_\mathcal{A})_\geq\ 
           \xrightarrow{\,\pi_\mathcal{A}\,}\ \Delta\,,
\]
where $\pi_\mathcal{A}$ is the linear projection of
Remark~\ref{R:moment_proj}, which restricts to give the moment map on
$(Y_\mathcal{A})_\geq$. 
Then linear precision of the components of ${\bf f}$ is simply the statement
that ${\bf f}$ is the inverse of the algebraic moment map.

\providecommand{\bysame}{\leavevmode\hbox to3em{\hrulefill}\thinspace}
\providecommand{\MR}{\relax\ifhmode\unskip\space\fi MR }
\providecommand{\MRhref}[2]{%
  \href{http://www.ams.org/mathscinet-getitem?mr=#1}{#2}
}
\providecommand{\href}[2]{#2}


\begin{thebibliography}{CLO97}

\bibitem[Ber75]{Be75}
D.~N. Bernstein, \emph{The number of roots of a system of equations},
  Funkcional. Anal. i Prilo\v zen. \textbf{9} (1975), no.~3, 1--4. \MR{55
  \#8034}

\bibitem[BGT97]{BGT97}
Winfried Bruns, Joseph Gubeladze, and Ng{\^o}~Vi{\^e}t Trung, \emph{Normal
  polytopes, triangulations, and {K}oszul algebras}, J. Reine Angew. Math.
  \textbf{485} (1997), 123--160. \MR{99c:52016}

\bibitem[Cox03]{Co02}
David Cox, \emph{What is a toric variety?}, 2003, Tutorial for Conference on
  Algebraic Geometry and Geometric Modeling, Vilnius, Lithuania, 29 July-2
  August.

\bibitem[CLO97]{CLO97}
David Cox, John Little, and Donal O'Shea, \emph{Ideals, varieties, and
  algorithms}, second ed., Undergraduate Texts in Mathematics, Springer-Verlag,
  New York, 1997, An introduction to computational algebraic geometry and
  commutative algebra. \MR{97h:13024}

\bibitem[BKK76]{BKK}
D.~Bernstein, A.~Kouchnirenko, and A.~Khovanskii, \emph{Newton polytopes}, Usp.
  Math. Nauk. \textbf{1} (1976), no.~3, 201--202, (in Russian).

\bibitem[Del03]{De03}
Claire Delaunay, \emph{Real structures on smooth compact toric surfaces}, 2003,
  Tutorial for Conference on Algebraic Geometry and Geometric Modeling,
  Vilnius, Lithuania, 29 July-2 August.

\bibitem[Ful93]{Fu93}
William Fulton, \emph{Introduction to toric varieties}, Annals of Mathematics
  Studies, vol. 131, Princeton University Press, Princeton, NJ, 1993, The
  William H. Roever Lectures in Geometry. \MR{94g:14028}

\bibitem[Mac2]{M2}
Daniel~R. Grayson and Michael~E. Stillman, \emph{Macaulay 2, a software system
  for research in algebraic geometry}, Available at
  {\tt http://www.math.uiuc.edu/Macaulay2/}.

\bibitem[SING]{SINGULAR}
G.-M. Greuel, G.~Pfister, and H.~Sch\"onemann, \emph{{\sc Singular} 2.0}, {A
  Computer Algebra System for Polynomial Computations}, Centre for Computer
  Algebra, University of Kaiserslautern, 2001, {\tt
 http://www.singular.uni-kl.de}.

\bibitem[Har92]{Ha92}
Joe Harris, \emph{Algebraic geometry}, Graduate Texts in Mathematics, vol. 133,
  Springer-Verlag, New York, 1992, A first course. \MR{93j:14001}

\bibitem[Koe93]{Ko93}
Robert~Jan Koelman, \emph{A criterion for the ideal of a projectively embedded
  toric surface to be generated by quadrics}, Beitr\"age Algebra Geom.
  \textbf{34} (1993), no.~1, 57--62. \MR{94h:14051}

\bibitem[Kra02]{Kr02}
Rimvydas Krasauskas, \emph{Toric surface patches}, Adv. Comput. Math.
  \textbf{17} (2002), no.~1-2, 89--133, Advances in geometrical algorithms and
  representations. \MR{2003f:65027}

\bibitem[Roj03]{Ro03}
J.~Maurice Rojas, \emph{Why polyhedra matter in non-linear equation solving},
  2003, Tutorial for Conference on Algebraic Geometry and Geometric Modeling,
  Vilnius, Lithuania, 29 July-2 August.

\bibitem[Stu96]{Sturmfels_GBCP}
Bernd Sturmfels, \emph{Gr\"obner bases and convex polytopes}, American
  Mathematical Society, Providence, RI, 1996. \MR{97b:13034}

\end{thebibliography}
\end{document}